\documentclass{amsart}
\usepackage{amsfonts}
\usepackage{amsmath,amscd}
\usepackage{amssymb}
\usepackage{amsthm}
\usepackage{newlfont}
\newcommand{\f}{\frac}

 \newtheorem{thm}{Theorem}
 
 \newtheorem{lem}[thm]{Lemma}
 
 \theoremstyle{definition}
 
 \theoremstyle{remark}

 \numberwithin{equation}{section}
\begin{document}

%
%\begin{center}Any suggestions are welcome !!!\end{center}
\title[Classifying $p$-groups via their multiplier]
 {Classifying $p$-groups via their multiplier}

\author[P. Niroomand]{Peyman Niroomand}

\email{niroomand@dubs.ac.ir}
\address{School of Mathematics and Computer Science\\
Damghan University of Basic Sciences, Damghan, Iran}


%\subjclass{}


\date{\today}

%\dedicatory{}

\begin{abstract} The author in $($On the order of Schur multiplier of non-abelian $p$-groups. J. Algebra (2009).322:
4479--4482$)$ showed that for any $p$-group $G$ of order $p^n$ there exists a nonnegative  integer $s(G)$
such that the order of Schur multiplier of $G$ is equal to
$p^{\f{1}{2}(n-1)(n-2)+1-s(G)}$. Furthermore, he characterized the
structure of  all non-abelian $p$-groups $G$ when $s(G)=0$. The
present paper is devoted to characterization  of all $p$-groups when
$s(G)=2$.
\end{abstract}

%%% ----------------------------------------------------------------------
\maketitle
%%% ----------------------------------------------------------------------

The concept of Schur multiplier, $\mathcal{M}(G)$, have been studied
by several authors, initiated by Schur in 1904. It is known that
the order of Schur multiplier of a given finite $p$-group of order
$p^n$ is equal to $p^{\f{1}{2}n(n-1)-t(G)}$ for some $t(G)\geq 0$ by
a result of Green \cite{gr}. It is of interest to know which
$p$-groups have the Schur multiplier of order
$p^{\f{1}{2}n(n-1)-t(G)}$, when $t(G)$ is in hand.

Historically, there are several papers trying to characterize the
structure of $G$ by just the order of its Schur multiplier.
In \cite{be} and \cite{zh}, Berkovich and Zhou classified the
structure of $G$ when $t(G)=0,1$ and $2$, respectively.

Later, Ellis in \cite{el} showed that having a new upper bound on
the order of  Schur multiplier of groups  reduces
characterization process of structure of $G$. He reformulated the
upper bound due to  Gasch\"utz at. al. \cite{ga} and classified with
a quite way to that of \cite{be, zh} the structure of $G$ when
$t(G)=3$.

 The result of \cite{ni} shows that  there exists a nonnegative integer $s(G)$ such
that $|\mathcal{M}(G)|=p^{\f{1}{2}(n-1)(n-2)+1-s(G)}$ which is a
reduction of Green's bound for any  given non-abelian $p$-group $G$
of order $p^n$. One can check that the structure of $G$ can be
characterized by using \cite[Main Theorem]{ni}, when
$t(G)=1,2,3$. Moreover, characterizing non-abelian $p$-groups by
 $s(G)$ can be significant since for instance the result of
\cite{ni} and \cite{ni1} emphasize that the number of groups with a
fixed $s(G)$ is more than that with fixed $s(G)$. Also the results
of \cite{ni2} and \cite{ni3} show  handling the $p$-groups
characterized by $s(G)=0,1$ may be caused to characterize the
structure of $G$ by $t(G)$.

In the present paper, we intend to classify the structure of all
 non-abelian $p$-groups when $s(G)=2$.

 Throughout this paper we use the following notations.
 \\$Q_8$: quaternion group of order $8$,
 \\$D_8$: dihedral group of order $8$,
 \\$E_1$: extra special $p$-group of order $p^3$ and exponent $p$,
 \\$E_2$: extra special $p$-group of order $p^3$ and exponent $p^2$ $(p\neq 2)$,
 \\${\mathbb{Z}}^{(m)}_{p^{n}}$: direct product of $m$ copies of the cyclic group of
order $p^n$,
\\$G^{ab}$: the abelianization of group $G$,
 \\$H\cdot K$: the central product of $H$ and $K$,
 \\$E(m)$: $E\cdot Z(E)$, where $E$ is an extra special $p$-group and $Z(E)$ is cyclic group of order $p^m$ $(m\geq 2)$,
 \\$\Phi(G)$: the Frattini subgroup of group $G$,
 \\Also, $G$  has the property $s(G)=2$ or briefly with $s(G)=2$ means  the order of its Schur multiplier is of order $p^{\f{1}{2}(n-1)(n-2)-1}$.

The following lemma is a consequence of \cite[Main Theorem]{ni}.
\begin{lem}\label{1} There is no $p$-group with $|G'|\geq p^3$ and $s(G)=2$.
\end{lem}
\begin{lem}\label{2} There is no $p$-group of order $p^n$ $(n\geq 5)$ when $G^{ab}$ is not elementary abelian and $s(G)=2$.
\end{lem}
\begin{proof} First suppose that $n=5$. By virtue of \cite[Theorem 3.6]{ni3}, the result follows.
In case $n\geq 6$, by invoking \cite[Lemma 2.3]{ni}, we have $|\mathcal{M}(G/G')|\leq p^{\f{1}{2}(n-2)(n-3)}$, and since $G/Z(G)$ is capable
the rest of proof is obtain by using \cite[Propsition 1]{el2}.
\end{proof}
\begin{lem}\label{3}Let $G$ be a $p$-group and $|G'|=p$ or $p^2$ with $s(G)=2$. Then $Z(G)$ is of exponent at most
$p^2$ and $p$, respectively.
\end{lem}
\begin{proof} Taking a cyclic central subgroup $K$ of order $p^k$ $(k\geq 3)$ and using \cite[Theorem 2.2]{j2}, we should have
\[\begin{array}{lcl}|\mathcal{M}(G)|\leq p^{-1}|G/K\otimes K|p^{\f{1}{2}(n-k)(n-k-1)}&\leq& p^{n-k-1}p^{\f{1}{2}(n-k)(n-k-1)}\vspace{.3cm}\\&\leq& p^{\f{1}{2}(n-1)(n-2)-2},\end{array}\]
which is a contradiction. In case $|G'|=p^2$, the result obtained
similarly.
\end{proof}
Lemma \ref{1} indicates when $G$ has the property $s(G)=2$, then
$|G'|\leq p^2$. First we suppose that $|G'|=p$.
\begin{thm}\label{4} Let $G$ be a $p$-group with
centre of order at most $p^2$ and $G^{ab}$ be elementary abelian of order $p^{n-1} $and $s(G)=2$. Then $G\cong E(2)$,
$E_2\times \mathbb{Z}_{p}$, $ Q_8$ or $H$, where $H$ is an extra
special $p$-group of order $p^{2m+1}$ $(m\geq 2)$.
\end{thm}
\begin{proof}First assume that $|Z(G)|=p$. Hence $G\cong Q_8$ or $G\cong H$, where $H$ is an
extra special $p$-group of order $p^{2m+1}$ $(m\geq 2)$ by a result of \cite[Theorem 3.3.6]{kar}.
 Now, assume that $|Z(G)|\geq p^2$.
Lemma \ref{3} and assumption show that
$Z(G)\cong\mathbb{Z}_{p}\times \mathbb{Z}_{p}$ or $
\mathbb{Z}_{p^2}$.

In case for which  $Z(G)$ is of exponent $p$, \cite[Lemma 2.1]{ni}
follows that $G\cong H\times \mathbb{Z}_{p}$. It is easily checked
that $H\cong E_2$ by using \cite[Theorems 2.2.10 and 3.3.6 ]{kar}.

In case $Z(G)$ is of exponent $p^2$, since $\Phi(G)=G'$,
\cite[Theorem 3.1]{j1} shows that
\[p^{\f{1}{2}(n-1)(n-2)}=|\mathcal{M}(G/\Phi(G))|\leq p~|\mathcal{M}(G)|,\] and hence $p^{\f{1}{2}(n-1)(n-2)-1}\leq|\mathcal{M}(G)|$.
On the other hand, Main Theorems of \cite{ni} and  \cite{ni1} imply that $|\mathcal{M}(G)|=p^{\f{1}{2}(n-1)(n-2)-1}$
since $Z(G)$ is cyclic of order $p^2$. Moreover, $G\cong E(2)$ by appealing \cite[Lemma 2.1]{ni}, as required.
\end{proof}

\begin{thm}\label{5}Let $G$ be a $p$-group, $G^{ab}$ be elementary abelian, $|G'|=p$ and $|Z(G)|\geq p^3$ be of exponent $p^2$. Then
$G=E(2)\times Z$, where $Z$ is an elementary abelian $p$-group.
\end{thm}
\begin{proof} It is known that $G=H\cdot Z(G)$ and $H\cap Z(G) = G'$ by virtue of \cite[Lemma 2.1]{ni}.
Now, for the  sake of clarity, we consider two cases.

Case1. First assume that $G'$ lies in a central subgroup $K$ of
exponent $p^2$. Therefore, one can check that there exists a central
subgroup $T$ such that $G=H\cdot K\times T\cong E(2)\times T$. Thus,
when $T$ is an elementary abelian $p$-group by using \cite[Theorem
2.2.10]{kar} and  Theorem \ref{4}, we have
\[\begin{array}{lcl}|\mathcal{M}(G)|&=&|\mathcal{M}(E(2))||\mathcal{M}(T)||E(2)^{ab}\otimes T|\vspace{.3cm}\\&=& 2m^2+m-1+
\f{1}{2}(n-2m-2)(n-2m-3)+2m+1(n-2m-2)\vspace{.3cm}\\&=&\f{1}{2}(n-1)(n-2)-1.\end{array}\]

In the case $T$ is not elementary abelian, a similar method and
\cite[Lemma 2.2]{ni} asserts that
\[|\mathcal{M}(G)|\leq p^{\f{1}{2}(n^2-5n+4)}\leq p^{\f{1}{2}(n-1)(n-2)-2}.\]

Case2.
 $G'$ has a complement $T$ in $Z(G)$, and  hence $G=H\times T$ where $T$ is not elementary abelian, and so by invoking
 \cite[Lemma 2.2]{ni} and \cite[Theorems 2.2.10 and 3.3.6]{kar}, $|\mathcal{M}(G)|\leq p^{\f{1}{2}(n-1)(n-2)-2}.$
 \end{proof}
 \begin{thm}\label{6}Let $G$ be a $p$-group of order $p^n$, $G^{ab}$ be elementary abelian of order $p^{n-1}$ and $Z(G)$ be
 of exponent $p$. Then $G$ has the property $s(G)=2$ if and only if it is isomorphic to one of the following groups.
 \[Q_8\times {\mathbb{Z}}^{(n-3)}_{2},  E_2\times{\mathbb{Z}}^{(n-3)}_{p}
 ~\text{or}~H\times{\mathbb{Z}}^{(n-2m-1)}_{p},
\]where $H$ is extra special of order $ p^{2m+1}$ $m\geq 2$.
 \end{thm}
 \begin{proof}It is obtained via Theorem \ref{4}, \cite[Theorems 2.2.10 and 3.3.6 ]{kar} and assumption.
 \end{proof}
 \begin{lem}Let $G$ be a $p$-groups of order $p^4$ and $|G'|=p$.
 Then $G$ has the property $s(G)=2$ if and only if $G$ is isomorphic
 to the one of the following groups.
 \begin{itemize}
 \item[(1)] $Q_8\times\mathbb{Z}_2$,
 \item[(2)] $\langle a,b~|~a^4=1, b^4=1,[a,b,a]=[a,b,b]=1,[a,b]=a^2b^2\rangle$
\item[(2)] $\langle a,b,c~|~a^2=b^2=c^2=1, abc=bca=cab\rangle$.
\item[(4)]$E_4\cong E_1(2)$,
\item[(5)]$E_2\times {\mathbb{Z}}_p$,
\item[(6)]$\langle a,b~|~a^{p^2}=1, b^p=1,[a,b,a]=[a,b,b]=1\rangle$,
\end{itemize}
 \end{lem}
 \begin{proof} It is obtained by using Theorems \ref{4},
 \ref{6} and a result of \cite[Lemma 3.5]{ni2}.
 \end{proof}
 The structure of all $p$-group of order $p^n$  is characterized with the property $s(G)=2$ and $|G'|=p$.
 Now, we may suppose that $|G'|=p^2$.
 \begin{lem}\label{7} There is no $p$-group of order $p^n$ $(n\geq 5)$ with $s(G)=2$, where
 $|G'|=p^2$ and $G'\nsubseteq Z(G)$.
 \end{lem}
 \begin{proof} First assume that $|Z(G)|=p^2$, since $Z(G)$ is
 elementary by Lemma \ref{4}. Let $K$ be a central subgroup of order
 $p$, such that $ |(G/K)'|=p^2$. It is seen that
\[|\mathcal{M}(G)|\leq |\mathcal{M}(G/K)||K\otimes G/(K\times
G')|\leq |\mathcal{M}(G/K)|~p^{n-3}\]
 by  \cite[Theorem 4.1]{j1}. On the other hand, \cite[Main
 Thoerems]{ni,ni1} imply that $ |\mathcal{M}(G/K)|\leq
 p^{\f{1}{2}(n-2)(n-3)-1}$, and hence
$|\mathcal{M}(G)|\leq p^{\f{1}{2}(n-1)(n-2)-2}$.

In case $|Z(G)|=p^3$, there exists a  central subgroup  $K$ of order
$p^2$ such that $G'\cap K=1$. The rest of proof is obtained similar
to the pervious case.

When $|Z(G)|=p$, since $G$ is nilpotent of class 3, the result is
deduced by \cite[Proposition 3.1.11]{kar}.
 \end{proof}
 \begin{thm}Let $G$ be a $p$-group of order $p^n$ $(n\geq 5)$ and $|G'|=p^2$ with $s(G)=2$. Then
 \[G\cong\mathbb{Z}_p\times \big({\mathbb{Z}}^{(4)}_p\rtimes_\theta
 \mathbb{Z}_p\big)~
 (p\neq2).\]
 \end{thm}
By the results of Lemmas \ref{4} and \ref{7}, we may assume that
$G'\subseteq Z(G)$ and $Z(G)$ is of exponent $p$. We consider three
cases relative to $|Z(G)|$.

Case1. Assume that $|Z(G)|=p^4$, there exists a central subgroup $K$
of order $p^2$ such that $K\cap G'=1$. \cite[Main Theorem]{ni}
implies that $|\mathcal{M}(G/K)|\leq p^{\f{1}{2}(n-3)(n-4)}$, and so
$|\mathcal{M}(G)|\leq p^{\f{1}{2}(n-1)(n-2)-1}$ due to \cite[Thorem
4.1]{j1}.

Case2. In the case $|Z(G)|=p^2$, we have $G'=Z(G)$. Moreover
\cite[Main Theorem]{ni1} deduces that $n\geq 6$ and so there exists
a central subgroup $K$ such that $G/K\cong H\times Z(G/K)$ where $H$
is a extra
 special $p$-groups of order $p^{2m+1}$ $m\geq 2$, thus
 \[|\mathcal{M}(G)|\leq p^{n-3}|\mathcal{M}(G/K)|\leq p^{n-3}p^{\f{1}{2}(n-1)(n-4)}\leq p^{\f{1}{2}(n-1)(n-2)-2}.\]
Case3. Now,  we may assume that $|Z(G)|=p^3$. Let $K$
be a complement of $G'$ in $Z(G)$, so  \cite[Main Theorem]{ni3}
asserts that $|\mathcal{M}(G/K)|\leq p^{\f{1}{2}(n-2)(n-3)}$. On the
other hand,
 \cite[Theorem 4.1]{j1} and assumption imply that
\[\begin{array}{lcl}p^{\f{1}{2}(n-1)(n-2)-1}=|\mathcal{M}(G)|&\leq& |\mathcal{M}(G/K)||K\otimes G/Z(G)|
\vspace{.3cm}\\&\leq& \big|\mathcal{M}(G/K)\big|p^{n-3},\end{array}\] so we
should have $|\mathcal{M}(G/K)|=p^{\f{1}{2}(n-2)(n-3)}$ and $G/Z(G)$
is elementary abelian. Now, since
$|\mathcal{M}(G/K)|=p^{\f{1}{2}(n-2)(n-3)}$ and $|(G/K)'|=p^2$, by
using \cite[Main Theorem]{ni1}, $G/K \cong
{\mathbb{Z}}^{(4)}_p\rtimes_\theta
 \mathbb{Z}_p~
 (p\neq2)$. Moreover, \cite[Proposition 1]{el2} and assumption show that $G^{ab}$ is elementary abelian. Hence, it is readily shown that
 \[G\cong\mathbb{Z}_p\times \big({\mathbb{Z}}^{(4)}_p\rtimes_\theta
 \mathbb{Z}_p\big)~
 (p\neq2).\]
\begin{thm}Let $G$ be a group of order $p^4$ with $s(G)=2$ and $|G'|=p^2$. Then $G$ is isomorphic to the one of the following groups.
\begin{itemize}
\item[(1)]$\langle
a,b~|~a^9=b^3=1,[a,b,a]=1,[a,b,b]=a^6,[a,b,b,b]=1\rangle$,
\item[(2)] $\langle a,b~|a^p=1, b^p=1,[a,b,a]=[a,b,b,a]=[a,b,b,b]=1\rangle (p\neq3)$.
\end{itemize}
\end{thm}
\begin{proof}The structure of these groups has been characterized in \cite[Lemma 3.6]{ni2}.
\end{proof}
We summarize all results as follows
\begin{thm}Let $G$ be a group of order $p^n$. Then $G$ has a
property $s(G)=2$ if and only if isomorphic to the one of the
following groups.
\end{thm}
\begin{itemize}
\item[(1)]$E(2)\times{\mathbb{Z}}^{(n-2m-2)}_{p}$,
\item[(2)]$E_2\times {\mathbb{Z}}^{(n-3)}_{p}$,
\item[(3)] $ Q_8\times {\mathbb{Z}}^{(n-3)}_{2}$
\item[(4)] $H\times{\mathbb{Z}}^{(n-2m-1)}_{p}$, where $H$ is an extra special
$p$-group of order $p^{2m+1}$ $(m\geq 2)$,
 \item[(5)] $\langle a,b~|~a^4=1, b^4=1,[a,b,a]=[a,b,b]=1,[a,b]=a^2b^2\rangle$
\item[(6)] $\langle a,b,c~|~a^2=b^2=c^2=1, abc=bca=cab\rangle$.
\item[(7)]$\langle a,b~|~a^{p^2}=1, b^p=1,[a,b,a]=[a,b,b]=1\rangle$,
\item[(8)]$\mathbb{Z}_p\times \big({\mathbb{Z}}^{(4)}_p\rtimes_\theta
 \mathbb{Z}_p\big)~
 (p\neq2)$,
\item[(9)]$\langle
a,b~|~a^9=b^3=1,[a,b,a]=1,[a,b,b]=a^6,[a,b,b,b]=1\rangle$,
\item[(10)] $\langle a,b~|a^p=1, b^p=1,[a,b,a]=[a,b,b,a]=[a,b,b,b]=1\rangle (p\neq3)$.
\end{itemize}

\end{document}